\documentclass[twoside,a4paper,11pt,limits]{amsart}
\usepackage{amsmath,amssymb}
\usepackage{paper_diening}
\usepackage{graphicx}
\usepackage{color}
\usepackage{latexsym}
\usepackage{amssymb}
\usepackage{amsfonts}
\usepackage{amsmath}
\newtheorem{theorem}{Theorem}[section]
\newtheorem{lemma}[theorem]{Lemma}
\newtheorem{corollary}[theorem]{Corollary}
\newtheorem{proposition}[theorem]{Proposition}
\newtheorem{remark}[theorem]{Remark}

\newtheorem{example}[theorem]{Example}

\numberwithin{equation}{section}
\newcommand{\meantmp}[2]{#1\langle{#2}#1\rangle}
\newcommand{\mean}[1]{\meantmp{}{#1}}

\newcommand{\R}{{\setR}}
\newcommand{\rn}{{\setR^n}}
\newcommand{\Rn}{{\setR^n}}
\newcommand{\RN}{{\setR^N}}
\newcommand{\RNn}{\setR^{n\times N}}

\newcommand{\diam}{\text{\rm diam}}

\newcommand{\dt}{\ensuremath{\,dt}}

\newcommand{\dx}{\ensuremath{\,dx}}
\newcommand{\dy}{\ensuremath{\,dy}}
\newcommand{\dz}{\ensuremath{\,dz}}

\newcommand{\comment}[1]{\vskip.3cm
\fbox{%
\parbox{0.93\linewidth}{\footnotesize #1}}
\vskip.3cm}

\newcommand{\tuomo}[1]{%
\textcolor[rgb]{1,0,0}{*}\textcolor[rgb]{0,1,0}{*}\textcolor[rgb]{0,0,1}{*}%
\textcolor[rgb]{1,0,0}{*}\textcolor[rgb]{0,1,0}{*}\textcolor[rgb]{0,0,1}{*}
\textcolor[rgb]{0.2,0.6,0.6}{Tuomo: #1}
\textcolor[rgb]{1,0,0}{*}\textcolor[rgb]{0,1,0}{*}\textcolor[rgb]{0,0,1}{*}%
\textcolor[rgb]{1,0,0}{*}\textcolor[rgb]{0,1,0}{*}\textcolor[rgb]{0,0,1}{*}}
\newcommand{\seb}[1]{{\textcolor[rgb]{0.00,0.00,1.00}{ #1}}}

\begin{document}

\title[Potential estimates for the $p$-Laplace system]{
Potential estimates for the $p$-Laplace system \\ with data in divergence form
%
}

\begin{abstract}{A pointwise  bound for local weak solutions to the
$p$-Laplace system is established in terms of data on the  right-hand side in divergence form. The relevant bound involves a Havin-Maz'ya-Wulff potential of the datum, and is a counterpart for data in divergence form of a classical result of \cite{KiMa}, that has  recently been extended  to systems in  \cite{KuMisys}. A local bound for oscillations is also provided.
%
These results allow for 
a unified
approach to regularity estimates for   broad classes of norms, including Banach function norms (e.g. Lebesgue, Lorentz and Orlicz norms), and norms depending on the oscillation of functions (e.g. H\"older, $\setBMO$ and, more generally, Campanato type norms). 
In particular,  new regularity properties are exhibited, and well-known
results  are easily recovered.   
%
%
%
}
\end{abstract}
\author{A.Cianchi and S.Schwarzacher}

\address{Andrea Cianchi, Dipartimento di Matematica e Informatica \lq\lq U.Dini", Universit\`{a} di Firenze, Viale Morgagni 67/A, 50134, Firenze, Italy}

\email{cianchi@unifi.it}

\address{Sebastian Schwarzacher, Department of mathematical analysis, Faculty of Mathematics and Physics,  Charles University in Prague,
Sokolovsk\'{a} 83, 186 75 Prague, Czech Republic}

\email{schwarz@karlin.mff.cuni.cz}

 \maketitle


\section{Introduction and main results}\label{intro}

The present paper deals with the regularity of local weak solutions to the  $p$-Laplace  system
\begin{align}
\label{eq:sysA}
  -\divergence(|\nabla \bfu|^{p-2} \nabla \bfu) &= -\divergence
  \bfF\, \qquad \hbox{in $\Omega$.}
\end{align}
 Here, $\Omega$ is an open set in $\mathbb R^n$, with $n \geq 2$, the exponent $p \in (1, \infty)$, the
function $\bfF \,:\, \Omega \to
\mathbb R^{N\times n}$, with $N
\geq 1$, is assigned, and $\bfu \,:\, \Omega \to \RN$ is the
unknown. 
\par We shall assume that $\bfF \in L^{p'}_{\rm loc}(\Omega)$ throughout, where $p'$ stands for the H\"older conjugate of $p$. This assumption guarantees that weak solutions to system \eqref{eq:sysA} are well defined.
Any weak solution $\bfu$ to system \eqref{eq:sysA} belongs, by definition, to the Sobolev space $W^{1,p}_{\rm loc}(\Omega)$. Basic  regularity properties of $\bfu$, such as membership in Lebsegue or H\"older spaces, according to whether $p \leq n$ or  $p>n$, can be immediately derived from this piece of information, via the standard Sobolev embedding theorem.  Additional regularity of $\bfF$, beyond local  $L^{p'}$-- integrability, is reflected into stronger regularity properties of any solution $\bfu$. 
\\ The regularity theory of solutions to $p$-Laplace type equations and systems has been the subject of a vast literature, starting from the second part of the last century. Classical fundamental contributions to this theory include \cite{Ura68}, \cite{Uhl77}, \cite{Iw83}, \cite{Ev}, \cite{Dibe}, \cite{Le}
\cite{Tolk}, \cite{ChenDiBe} \cite{DiBMan93}.
\par Instead of focusing on estimates for specific norms, here we offer a precise  pointwise estimate for any weak solution $\bfu$ to \eqref{eq:sysA}, and an estimate for its oscillations in integral form.
They provide us with  versatile tools for the proof of norm bounds in a wide range of function spaces, as shown in the last part of the paper.  

Our first  main result is  contained in Theorem \ref{mainpointwise}, and amounts to a  pointwise bound for $\bfu$ in terms of a Havin-Maz'ya-Wulff potential of $\bfF$. This can be regarded as a version, for right-hand sides in divergence form, of the pointwise bound established in \cite{KiMa} for equations whose right-hand side is a function (or, more generally, a measure), and of its recent extension to systems from \cite{KuMisys}. The latter paper also contains a parallel pointwise estimate for the gradient via a Riesz potential of the datum, which carries over to systems a result from \cite{KuMiARMA}. Earlier contributions along a similar line of research are \cite{Mingione2011} and \cite{DuzMing2011}. Pointwise gradient estimates for solutions to systems with right-hand side in the form of \eqref{eq:sysA} are proved in \cite{BCDKS1} -- see also \cite{BCDKS2}.  Let us add that  rearrangement bounds for solutions to boundary value problems for $p$-Laplace type elliptic equations with the same kind of right-hand side can be found in \cite{BFM}. 
Pointwise  estimates, in rearrangement form,  for
the gradient  are in  \cite{AFT}, \cite{ACMM} and  \cite{CiaMaz}.
\par
Recall that, given $s
>1$ and $\sigma
>0$,   the truncated Havin-Maz'ya-Wulff  potential 
$\bfW_{\alpha,s}^R\bff$ of an $\mathbb R^m$-valued function $\bff \in L^1_{\rm loc}(\Omega )$, with $m \geq 1$, is
defined as
\begin{align}\label{wulff}
 \bfW_{\alpha,s}^R\bff (x)=\int_0^R\bigg(r^{\alpha s}\dashint_{B_r(x)}\abs{\bff}\dy\bigg)^\frac1{s-1}\frac{d\,r}{r}
\end{align}
for every $x \in \Omega$ and $R>0$ such that $B_R(x) \subset 
\Omega$. Here, $B_r(x)$ denotes the ball centered at $x$, with radius $r$, and $\dashint_E$ stands for the averaged integral $\frac 1{|E|}\int _E$, where $|E|$ denotes the Lebesgue measure of a set $E\subset \rn$.

%


\begin{theorem}
\label{mainpointwise} Let $n \geq 2$, $N \geq 1$ and $p \in (1, \infty)$. Assume that
$\Omega $ is an open set in $\Rn$, and that $\bfF \in L^{p'}_{\rm
loc}(\Omega)$. Let $\bfu \in W^{1,p}_{\rm loc}(\Omega)$ be a
local weak solution to system \eqref{eq:sysA}.
There exists a
constant  $C=C(n, N, p)$ such
that
\begin{align}\label{main1}
 \abs{\bfu(x)}
&\leq C 
\bfW_{\frac{p}{p+1},p+1}^R\big(\abs{\bfF}^{p'}\big)(x)+C \dashint_{B_{R}(x)}\abs{\bfu}\dy
\end{align}
for a.e. $x \in \Omega$ and  every $R>0$ such that $B_{R}(x)\subset
{\Omega}$.    Moreover, a point $x \in \Omega$
 is a Lebesgue point of
$\bfu$ whenever the right-hand side
of \eqref{main1} is finite for some $R>0$.
%
%
%
%
%
\end{theorem}

\begin{remark}
{\rm Assume that $\Omega = \rn$, and that $\bfu$ is a weak solution to system \eqref{eq:sysA}  decaying so fast  at infinity that 
\begin{equation}\label{lim0}
\lim _{R \to \infty} \dashint_{B_{R}(x)}\abs{\bfu}\dy =0\,.
\end{equation}
Then, passing to the limit as $R \to \infty$ in inequality \eqref{main1} yields
%
%
\begin{equation}\label{mainrn}
\abs{\bfu(x)}
\leq C  
\bfW_{\frac{p}{p+1},p+1}\big(\abs{\bfF}^{p'}\big)(x) \quad \hbox{for a.e. $x \in \rn$.}
\end{equation} 
Here, $\bfW_{\alpha,s}\bff$ denotes the potential of a function $\bff \in L^1_{\rm loc}(\rn )$ given by the integral on the right-hand side 
of \eqref{wulff}, with $R$ replaced by $\infty$. Namely,
\begin{align}\label{wulffinf}
 \bfW_{\alpha,s}\bff (x)=\int_0^\infty\bigg(r^{\alpha s}\dashint_{B_r(x)}\abs{\bff}\dy\bigg)^\frac1{s-1}\frac{dr}{r}
\end{align}
for  $x \in \rn$. Note that condition \eqref{lim0} is fulfilled, for instance, if $p \in (1,n)$, 
$\nabla \bfu \in L^p(\rn)$ and $|\{x\in \rn: |\bfu(x)|>t\}|<\infty$ for every $t>0$. Indeed, \eqref{lim0} follows via
the H\"older and the classical Sobolev inequality, which holds 
under these assumptions on $\bfu$.
}
\end{remark}



\begin{remark} {\rm
In the case when $p=2$, system \eqref{eq:sysA} reduces to the classical Poisson system
\begin{equation}\label{poisson}
- \Delta \bfu = - \divergence \bfF\,,
\end{equation}
and 
inequality
\eqref{main1}  reads
\begin{equation}\label{linear}
\abs{\bfu(x)}\leq
c\int_0^R\bigg(\dashint_{B_r(x)}\abs{\bfF}^2\,dy\bigg)^\frac12\,dr+c\dashint_{B_{R}(x)}\abs{\bfu}\dy.
\end{equation}
The operator acting on $|\bfF|$ in \eqref{linear} bounds, via
H\"older's inequality,  the operator $\int
_0^R\dashint_{B_r(x)}\abs{\bfF}\,dy\, dr$, which is, in turn, 
equivalent to a truncated Riesz potential of order $1$ of $|\bfF|$.
In fact, if $\Omega = \rn$, and $F$ is regular enough, then the solution $\bfu$ to \eqref{poisson}, that decays to $0$ near infinity, admits the representation formula 
 $$\bfu (x) = c \int _{\rn} \bfF (y) 
\frac{x-y}{|x-y|^n}\, dy \qquad \hbox{for $x \in \rn$,}
$$
\color{black}
for a suitable constant $c=c(n)$.
The reason why Theorem \ref{mainpointwise} cannot yield inequality \eqref{linear}, with 
$\int_0^R\big(\dashint_{B_r(x)}\abs{\bfF}^2\,dz\big)^\frac12\,dr$
replaced by $\int
_0^R\dashint_{B_r(x)}\abs{\bfF}\,dy\, dr$, is that
 data $\bfF$, that do not belong to
$L^{p'}_{\rm loc}(\Omega)$, are not included in our analysis, which is nonlinear in nature.  
This is
related to a gap between the linear Calder\'on-Zygmund theory for the Laplacian, and its  nonlinear counterpart for the $p$-Laplacian, and hence to a well known open problem in the latter about right-hand sides $\bfF \in L^q_{\rm loc}(\Omega)$ with $q<p'$. See \cite{Iwa92},  \cite{Le93}  and \cite{KiLe} in this connection, and also the recent contributions
\cite{BulDieSch16,BulSch16}.
}

\end{remark}

\par\noindent

\color{black}

\smallskip
Theorem \ref{mainpointwise} applies, in principle, to regularity estimates for solutions to \eqref{eq:sysA} in any norm
 depending only on the size of functions, or, more precisely, in any norm that is monotone under poitwise domination of functions.
 They are usually called Banach function norms in the literature \cite{BS}. This class of norms includes various customary instances, such as the (possibly weighted) Lebesgue norms, the Orlicz norms and the Lorentz norms.
\par Bounds in norms depending on oscillations of functions, such as H\"older, BMO and,  more generally, Campanato type norms can be derived from the next result. Its content is an estimate for the oscillation of solutions to \eqref{eq:sysA} on balls. The latter   is also our first step in the proof of Theorem \ref{mainpointwise}. In the statement, and in what follows,  the notation $\mean{\bff}_E$ is  used, when convenient, to denote the integral average $\dashint_E \bff(x)\, dx$ of a function $\bff$ over a set $E$.

\begin{theorem}
\label{mainosc} Let $n \geq 2$, $N \geq 1$ and $p \in (1, \infty)$. Assume that
$\Omega $ is an open set in $\Rn$, and that $\bfF \in L^{p'}_{\rm
loc}(\Omega)$. Let $\bfu \in W^{1,p}_{\rm loc}(\Omega)$ be a
local weak solution to system \eqref{eq:sysA}.
There exists a
constant $C=C(n, N, p)$ such that 
\begin{align}\label{main2}
 \dashint_{B_r(x)}\abs{\bfu(y)-\mean{\bfu}_{B_r (x)}} dy
&\leq C r \bigg(\int_r^{R}\bigg( \dashint_{B_\rho (x)}
\abs{\bfF-\mean{\bfF}_{B_\rho (x)}}^{p'} \dy
\bigg)^\frac{1}{p'}\frac{d\rho}{\rho}\bigg)^{\frac 1{p-1}} +  C r
\dashint_{B_{R}(x)}\abs{ \nabla\bfu}\dy
\end{align}
for a.e. $x \in \Omega$,  every $R>0$ such that $B_{R}(x)\subset
{\Omega}$ and every $r\in (0,R]$.
%
%
%
%
\end{theorem}

\begin{remark}\label{pointnew}{\rm A close inspection of the proof of Theorem \ref{mainpointwise} will reveal that inequality \eqref{main1} can be slightly improved, in that the Wulff potential 
$
\bfW_{\frac{p}{p+1},p+1}^R\big(\abs{\bfF}^{p'}\big)$ on its right-hand side can be replaced with a smaller nonstandard potential, introduced in \cite{BCDKS1}. The relevant potential involves the oscillation of $\bfF$ on balls, as in \eqref{main2}, instead of just its $L^{p'}$ averages, as in $\bfW_{\frac{p}{p+1},p+1}^R\big(\abs{\bfF}^{p'}\big)$. The resulting inequality reads 
\begin{align}\label{pointosc}
\abs{\bfu(x)} \leq 
C\int_0^R\bigg( \dashint_{B_{\rho}(x)}
\abs{\bfF-\mean{\bfF}_{B_\rho (x)}}^{p'} \dy
\bigg)^\frac{1}{p}d\rho + C \dashint_{B_{R}(x)}\abs{
\bfu}\dy
\end{align}
for a.e. $x \in \Omega$ and  every $R>0$ such that $B_{R}(x)\subset
{\Omega}$.  Clearly, inequality \eqref{pointosc} implies \eqref{main1}, 
since dropping the expression $\mean{\bfF}_{B_\rho (x)}$ in the first integral of its right-hand side results in  $\bfW_{\frac{p}{p+1},p+1}^R\big(\abs{\bfF}^{p'}\big)(x)$.
}
\end{remark}
 
 The proofs of Theorems \ref{mainpointwise} and \ref{mainosc} are accomplished in the next section. Tools playing a role in our approach include  a recent pointwise estimate for
a sharp maximal function of the gradient of solutions to system \eqref{eq:sysA} from \cite{BCDKS1},  suitable versions of  classical results in the regularity theory of the $p$-Laplacian  system, as well as   certain weighted inequalities of Hardy type for functions of one variable.  In Section 3, Theorem \ref{mainpointwise} is exploited, in combination with apropos results on the boundedness of nonlinear potentials, to establish regularity bounds for solutions in Lorentz and Orlicz spaces. Bounds in Campanato type spaces, and ensuing estimates in spaces of uniformly continuous functions are derived in the same section, via Theorem  \ref{mainosc}.

\section{Proofs of Theorems \ref{mainpointwise} and \ref{mainosc}}\label{proofs}

The  Hardy type inequalities contained in the next lemma
will be needed in our proof of Theorem~\ref{mainpointwise}. In what follows,
a function $\varphi : (0, L ) \to [0, \infty)$, with $L \in (0,
\infty]$ will be called quasi-increasing if there exists a constant
$k \geq 1$ such that
\begin{equation}\label{quasi}\varphi (r) \leq k \varphi (s) \qquad \hbox{if $0< r \leq s<L$.}
\end{equation}
Note that, if the function $\varphi$ fulfills inequality \eqref{quasi}, then the function $\psi : (0, L ) \to [0, \infty)$ associated
with $\varphi$ as
$$  \psi (s) = \sup _{0<r<s} \varphi (r) \quad \hbox{for $s \in (0, L)$,}$$
 is  non-decreasing, and  satisfies
the inequalities
\begin{equation}\label{quasi1} \varphi (s) \leq \psi (s) \leq k \varphi
(s) \quad \hbox{for $s \in (0, L)$.}
\end{equation}

\medskip

 \begin{lemma}\label{hardy} Let $\alpha \in \mathbb R$.
\par\noindent
{\rm (i)} Assume that $q \in [1,\infty)$. Then there exists a constant
$C=C(q, \alpha)$ such that
\begin{equation}\label{hardy1}
 \bigg(\int _0^\infty \bigg(\int _s^\infty \varphi (r) r^\alpha \,dr
 \bigg)^q\, ds\bigg)^{\frac 1q} \leq C \bigg(\int _0^\infty \varphi (s)^q
 s^{q(\alpha +1)}\, ds \bigg)^{\frac 1q}
 \end{equation}
 for every measurable function $\varphi : (0, \infty ) \to [0, \infty)$.
 \par\noindent
{\rm (ii)} Assume that $q \in (0,1)$. If $\alpha < -1 - \frac 1q$, then
there exists a constant $C=C(q, \alpha, k)$ such that inequality
\eqref{hardy1} holds
 for every quasi-increasing function $\varphi : (0, \infty ) \to [0, \infty)$ with constant $k$ as in  \eqref{quasi}.
 If $-1 - \frac 1q \leq \alpha < -1 $ and $0<a< \infty$, then
 there exists a constant $C=C(q, \alpha, k)$ such that
\begin{equation}\label{hardy2}
 \bigg(\int _0^a \bigg(\int _s^a \varphi (r) r^\alpha \,dr
 \bigg)^q\, ds\bigg)^{\frac 1q} \leq C \bigg(\int _0^{2a} \varphi (s)^q
 s^{q(\alpha +1)}\, ds \bigg)^{\frac 1q}
 \end{equation}
 for every quasi-increasing function $\varphi : (0, 2a ) \to [0, \infty)$ with constant $k$ as in  \eqref{quasi}.
\color{black}
\end{lemma}

\par\noindent
\textbf{Proof}. \emph{Part (i)}. This is a classical Hardy inequality -- see e.g.
\cite[Theorem 1.3.2/3]{Mazbook}.
\par\noindent
\emph{Part (ii)}. By inequalities \eqref{quasi1}, it   suffices
to prove the statement for non-decreasing
  functions $\varphi$.
\\ Let $0<a \leq b \leq \infty$. A characterization of weighted Hardy type inequalities
for monotone functions tells us that the inequality 
\begin{equation}\label{hardy4}
 \bigg(\int _0^a  \bigg(\int _s^a \varphi (r)   r^\alpha \,dr
 \bigg)^q\, ds\bigg)^{\frac 1q} \leq C \bigg(\int _0^{b} \varphi (s)^q
 s^{q(\alpha +1)}\, ds \bigg)^{\frac 1q}
 \end{equation}
holds  for every non-decreasing function $\varphi : (0, b) \to [0,
 \infty)$ if and only if
\begin{equation}\label{hardy3}
 \bigg(\int _0^a   \bigg(\int _{\max \{t,s\}}^b  \chi_{(0,a)}(r) r^\alpha \,dr
 \bigg)^q\, ds\bigg)^{\frac 1q} \leq C \bigg(\int _t^{b}
 s^{q(\alpha +1)}\, ds \bigg)^{\frac 1q} \quad \hbox{for every $t \in (0,
 b)$,}
 \end{equation}
with the same constant $C$ (see e.g. \cite[Theorem
3.3]{HeinigMaligranda}).  If $\alpha < -1 - \frac 1q$, then 
\begin{align}\label{hardy10}
\frac{\big(\int _0^\infty  \big(\int _{\max \{t,s\}}^\infty   r^\alpha \,dr
 \big)^q\, ds\big)^{\frac 1q}
}{\big(\int _t^{\infty}
 s^{q(\alpha +1)}\, ds \big)^{\frac 1q}
}
=C \qquad \hbox{for $t >0$,}
\end{align}
for some constant $C=C(\alpha , q)$. Thus, inequality \eqref{hardy3}, and hence also inequality \eqref{hardy4}, follows with $a=b=\infty$.
\\ Suppose next that $\alpha = -1 - \frac 1q$ and $b=2a< \infty$. Straightforward computations show that  there exist constants $C$ and $C'$, depending on $q$, such that
\begin{align}\label{hardy11}
\frac{\big(\int _0^a  \big(\int _{\max \{t,s\}}^a  \chi_{(0,a)}(r) r^\alpha \,dr
 \big)^q\, ds\big)^{\frac 1q}
}{\big(\int _t^{2a}
 s^{q(\alpha +1)}\, ds \big)^{\frac 1q}
}
& = C \frac{\big[t(t^{-\frac 1q}-a^{-\frac 1q})^q + \int_{t/a}^1(s^{-\frac 1q}-1)^q\, ds\big]^{\frac 1q}}{\big(\log \frac {2a}t\big)^{\frac 1q}
}
\\ \nonumber & \leq 
C \frac{\big[1 + \int_{t/a}^1 \frac{ds}s\big]^{\frac 1q}}{\big(\log \frac {2a}t\big)^{\frac 1q}
}
=C \frac{\big[1 + \log \frac at\big]^{\frac 1q}}{\big(\log \frac {2a}t\big)^{\frac 1q}
}
\\ \nonumber &\leq C' \Bigg[\frac{1}{\big(\log \frac {2a}t\big)^{\frac 1q}}+ \bigg(\frac{\log \frac at}{\log \frac {2a}t}\bigg)^{\frac 1q}\Bigg]\leq C' \Bigg[\frac{1}{(\log 2)^{\frac 1q}} + 1\Bigg]
\end{align}
for $t \in (0, a)$. Since the leftmost side of the chain \eqref{hardy11} vanishes for $t \in [a, 2a)$, inequality \eqref{hardy3}, and hence also inequality \eqref{hardy4}, follows with $b=2a$.
\\
 Finally, assume that $-1 - \frac 1q < \alpha < -1$ and $b=2a< \infty$. Then, there exist  constants $C=C(\alpha , q)$ and $C'=C'(\alpha , q)$ such that
 \begin{align}\label{hardy12}
\frac{\big(\int _0^a  \big(\int _{\max \{t,s\}}^a  \chi_{(0,a)}(r) r^\alpha \,dr
 \big)^q\, ds\big)^{\frac 1q}
}{\big(\int _t^{2a}
 s^{q(\alpha +1)}\, ds \big)^{\frac 1q}
}
 & = C a^{\alpha +1 + \frac 1q}\frac{\big[(\frac ta)^{q(\alpha +1) +1}(1 - (\frac at)^{\alpha +1})^q + \int _{ t/a}^1(s^{\alpha +1} -1)^q\, ds\big]^{\frac 1q}}{[(2a)^{q(\alpha +1) +1} - t^{q(\alpha +1) +1}]^{\frac 1q}}
 \\ \nonumber & \leq 
 C'
 a^{\alpha +1 + \frac 1q}\frac{\big[(\frac ta)^{q(\alpha +1) +1}(1 - (\frac at)^{\alpha +1})^q + (1- (\frac ta)^{q(\alpha +1) +1})
 \big]^{\frac 1q}}{[(2a)^{q(\alpha +1) +1} - t^{q(\alpha +1) +1}]^{\frac 1q}}
 \\ \nonumber & \leq 
 C'
 \frac{a^{\alpha +1 + \frac 1q} 2^{\frac 1q}}{a^{\alpha +1 +\frac 1q}} = C' 2^{\frac 1q}
 \end{align}
 for $t \in (0, a)$. The leftmost side of  \eqref{hardy12} vanishes for $t \in [a, 2a)$, and therefore inequalities \eqref{hardy3} and  \eqref{hardy4}, hold with $b=2a$.
\color{black}
\qed

\bigskip\par
In what follows, we shall make repeatedly   use of the inequality
\begin{equation}\label{mv}
\int _E |\bff - \mean{\bff}_E|\, \dx \leq 2 \int _E |\bff - \bfc|\, \dx
\end{equation}
for every measurable set $E \subset \rn$, every function $\bff$ in $E$, and every vector $\bfc \in \rn$. In particular, it plays a role in the proof of the next lemma, whose objective is an estimate for the difference between the mean values of a function over balls. In its statement and proof, all balls are concentric. Since the center is irrelevant, it will be dropped in the notation.    

\begin{lemma}
\label{lem:iter}
Let $R>0$ and let $\bff\in L^{{1}}(B_R)$.  Then 
\begin{equation}\label{tele1}
\abs{\mean{\bff}_{ B_r}-\mean{\bff}_{B_R}}
\leq 2^{2n+2}\int^R_r \dashint_{B_\rho}\abs{\bff-\mean{\bff}_{ B_\rho}}\, dx\frac{d\rho}\rho \qquad \hbox{for  $r\in(0,R]$}, 
\end{equation}
and
\begin{equation}\label{tele2}
\abs{\mean{\abs{\bff}}_{ B_r}-\mean{\abs{\bff}}_{B_R}}
\leq 2^{2n+3}\int^R_r \dashint_{B_\rho}\abs{\bff-\mean{\bff}_{ B_\rho}}\, dx\frac{d\rho}\rho \qquad \hbox{for  $r\in(0,R]$.}
\end{equation}
\end{lemma}
\par\noindent {\bf Proof}.
Since inequalities \eqref{tele1} and \eqref{tele2} are scale-invariant, we may assume, without  loss of generality, that $R=1$.
To begin with, observe that
\begin{align*}
  \abs{\mean{\bff}_{B_{\frac 12}}-\mean{\bff}_{B_1}}\leq
  \dashint_{B_{\frac 12}}\abs{\bff-\mean{\bff}_{B_1}}dx \leq  2^n\dashint_{B_1}\abs{\bff-\mean{\bff}_{B_1}}\,dx\,.
\end{align*}
Applying this inequality with $B_1$ replaced by $B_{2^{-i}}$, for $i=0, \dots , m-1$, and adding the resulting inequalities yield
\begin{align}
  \label{eq:iteration}
  \abs{\mean{\bff}_{B_{2^{-m}}}-\mean{\bff}_{B_1}}\leq 2^n\sum_{i=0}^{m-1}
    \dashint_{B_{2^{-i}}}\abs{\bff-\mean{\bff}_{B_{2^{-i}}}}\, dx.
\end{align}
Let $\theta\in (2^{-m},2^{1-m}]$. Owing to inequality \eqref{mv}, 
\begin{equation}\label{tele3}
\dashint_{B_{2^{-m}}}\abs{\bff-\mean{\bff}_{B_{2^{-m}}}}\,dx\leq 2^{n+1}\dashint_{  B_\theta }\abs{\bff - \mean{\bff}_{B_\theta \color{black}}}\,dx\leq 2^{2n+2}\dashint_{ B_{2^{1-m}}}\abs{\bff - \mean{\bff}_{B_{2^{1-m}}}}\,dx.
\end{equation}
On exploiting inequality \eqref{tele3}, one can show that 
\begin{equation}\label{tele4}
\dashint_{B_{2^{-m}}}\abs{\bff -\mean{\bff}_{B_{2^{-m}}}}\,dx\leq 2^{n+1}\dashint_{2^{-m}}^{2^{1-m}}\dashint_{  B_\rho }\abs{\bff -\mean{\bff}_{ B_\rho }}\,d  x\,d\rho\leq 2^{n+ 2 \color{black}}\int_{2^{-m}}^{2^{1-m}}\dashint_{  B_\rho }\abs{\bff-\mean{\bff}_{B_\rho }}\,dx\frac{d\rho}{\rho}.
\end{equation}
Now, given $r \in (0, 1]$, let $m \in \mathbb N$ be such that   $r\in (2^{-m},2^{1-m}]$. Thanks to 
inequality \eqref{tele4}, we have that
\begin{align*}
\abs{\mean{\bff}_{  B_r}-\mean{\bff}_{B_1}}&\leq\abs{\mean{\bff}_{ B_r}-\mean{\bff}_{B_{2^{1-m}}}}+\abs{\mean{\bff}_{  B_{2^{1-m}}}-\mean{\bff}_{B_1}}
\\
&\leq 2^n\dashint_{B_{2^{1-m}}}\abs{\bff-\mean{\bff}_{B_{2^{1-m}}}}\dx+2^n\sum_{i=0}^{m-2}   \dashint_{B_{2^{-i}}}\abs{\bff -\mean{\bff}_{B_{2^{-i}}}}\dx
\\
&\leq 2^n\sum_{i=0}^{m-1}   \dashint_{B_{2^{-i}}}\abs{\bff -\mean{\bff}_{B_{2^{-i}}}}\dx
\leq  2^{2n+{2}}\int_{2^{1-m}}^1\dashint_{  B_\rho}\abs{\bff -\mean{\bff}_{ B_\rho}}\dx\frac{d\rho}{\rho}
\\
&\leq 2^{2n+ 2 \color{black}}\int_{r}^1\dashint_{  B_\rho }\abs{\bff -\mean{\bff}_{ B_\rho}}\dx\frac{d\rho}{\rho}.
\end{align*}
This concludes the proof of inequality \eqref{tele1}. Inequality \eqref{tele2} follows via the same argument, combined with the fact that, by inequality \eqref{mv}, 
\begin{align}
\label{eq:abs}
 \dashint_{B_1}\abs{\abs{\bff}-\mean{\abs{\bff}_{B_1}}}\,dx\leq 2\dashint_{B_1}\abs{\abs{\bff}-\abs{\mean{\bff}_{B_1}}}\,dx\leq 2\dashint_B\abs{\bff-\mean{\bff}_{B_1}}\,dx. \qquad \qquad \qquad \qquad \qed
\end{align}

\smallskip
\par 
We are now in a position to prove our main results.

\smallskip
\par\noindent
\textbf{Proof of Theorem~\ref{mainosc}}.  Let $\bfu$ be a local weak solution to system \eqref{eq:sysA}. 
This means that
$\bfu \in W^{1,p}_{\loc}(\Omega)$
and
\begin{equation}\label{weaksol}
\int _{\Omega '} |\nabla \bfu|^{p-2}\nabla \bfu \cdot \nabla \bfphi
\, dx = \int _{\Omega '} \bfF \cdot \nabla \bfphi\,dx
\end{equation}
for every function $\bfphi \in W^{1,p}_0(\Omega ' )$, and every open
set $\Omega ' \subset \subset \Omega$. Here, the dot $\lq\lq \cdot
"$ stands for scalar product.
Let us set 
$$\bfA(\nabla\bfu) = |\nabla \bfu |^{p-2} \nabla \bfu.$$  Fix any $x \in \Omega$ and $R>0$ such that $B_{R}(x) \subset \Omega$. Let $r \in (0, R/2]$.
 By a standard \Poincare inequality on balls, and  H\"older's inequality,
\begin{align}
\label{ponc}
\bigg(\dashint_{B_r(x)}\abs{\bfu(y)-\mean{\bfu}_{B_r}(x)} \,dy\bigg)^{p-1}&\lesssim
\bigg(r \dashint_{B_r(x)}\abs{\nabla\bfu(y)}\, dy\bigg)^{p-1} \\
& \nonumber \lesssim r^{p-1}\bigg(\dashint_{B_r(x)}\abs{\bfA(\nabla\bfu(y))}^{\min\set{p',2}}\, dy\bigg)^\frac{1}{\min\set{p',2}}.
\end{align}
Here, and in what follows, the relation $\lq\lq \lesssim "$  between two expressions means that the former is bounded by the latter, up to multiplicative constants independent of the relevant variables involved. 
Inequality \cite[(3.16)]{BCDKS1} tells us that
\begin{multline}\label{june46}
\sum_{i=0}^{k}  \bigg(\dashint_{B_{\theta^i {\frac R2}}(x)}\abs{\bfA(\nabla
\bfu) -\mean{\bfA(\nabla \bfu)}_{B_{\theta^i {\frac R2}}(x)
}}^{\min\set{p',2}}\dy \bigg)^\frac1{\min\set{p',2}}\\ \lesssim
\bigg(\dashint_{B_{\frac R2}(x)}\abs{\bfA(\nabla \bfu)-\mean{\bfA(\nabla
\bfu)}_{ B_{\frac R2}(x) }}^{\min\set{p',2}}\dy\bigg)^\frac1{\min\set{p',2}}
+  \sum_{i=0}^{k-1}\bigg(\dashint_{B_{\theta^i {\frac R2}}(x)}
\abs{\bfF-\mean{\bfF}_{B_{\theta^i
{\frac R2}}(x)}}^{p'}\dy\bigg)^\frac{1}{p'}
\end{multline}
for every $\theta\in (0,1)$ and $k\in \setN$.
 Now, choose $k \in \setN$, such
that $\theta^{k+1} R\leq r\leq \theta^{k} R$ in \eqref{june46}. Via a telescope sum
argument, one can infer that
\begin{align}\label{2017-1}
&\bigg(\dashint_{B_r(x)}\abs{\bfA(\nabla\bfu)}^{\min\set{p',2}} dy\bigg)^\frac{1}{\min\set{p',2}}
\lesssim
\bigg(\dashint_{ B_{\theta^{k}{\frac R2}}(x)}\abs{\bfA(\nabla \bfu)}^{\min\set{p',2}}\dy \bigg)^\frac{1}{\min\set{p',2}}\\ \nonumber 
&\quad \lesssim
\bigg(\dashint_{ B_{\theta^{k}{\frac R2}}(x)}\abs{\bfA(\nabla \bfu)-\mean{\bfA(\nabla \bfu)}_{ B_{\theta^{k}{\frac R2}}(x)}}^{\min\set{p',2}}\dy \bigg)^\frac{1}{\min\set{p',2}}
\\ \nonumber 
&\qquad+ \sum_{i=1}^k\abs{\mean{\bfA(\nabla \bfu)}_{ B_{\theta^{i-1}{\frac R2}}(x)}-\mean{\bfA(\nabla \bfu)}_{ B_{\theta^{i}{\frac R2}}(x)}}
+ \abs{\mean{\bfA(\nabla \bfu)}_{ B_{\frac R2}(x)}}\\ \nonumber
&\quad\lesssim \sum_{i=0}^k\bigg(\dashint_{ B_{\theta^{i}{\frac R2}}(x)}\abs{\bfA(\nabla \bfu)-\mean{\bfA(\nabla \bfu)}_{ B_{\theta^{i}{\frac R2}}(x)}}^{\min\set{p',2}}\dy \bigg)^\frac{1}{\min\set{p',2}} + \abs{\mean{\bfA(\nabla \bfu)}_{ B_{\frac R2}(x)}}\\ \nonumber 
&\quad \lesssim \sum_{i=0}^{k-1}\bigg(\dashint_{B_{\theta^i {\frac R2}}(x)}
\abs{\bfF-\mean{\bfF}_{B_{\theta^i
{\frac R2}}(x)}}^{p'}\dy\bigg)^\frac{1}{p'} +\bigg(\dashint_{ B_{\frac R2}(x)}\abs{\bfA(\nabla \bfu)}^{\min\set{p',2}}\dy \bigg)^\frac{1}{\min\set{p',2}},
\end{align}
where the second inequality holds by an iterated use of the triangle inequality, the third one by H\"older's inequality, and the last one by \eqref{june46}. A standard 
reverse H\"older's inequality ensures that
\begin{align}
\label{revH}
\bigg(\dashint_{B_{\frac R2}(x)}\abs{\nabla \bfu}^p\dx\bigg)^\frac1p\lesssim \dashint_{B_{R}(x)}\abs{\nabla \bfu}\dx+\bigg(\dashint_{B_{R}(x)}\abs{\bfF-\mean{\bfF}_{B_{R}(x)}}^{p'}\dx\bigg)^\frac1{p},
\end{align}
see e.g.
~\cite[Lemma 3.2]{DuzMin10a}, or \cite[Lemma 3.2, Lemma 3.3.]{DieKapSch11}. 
An application of H\"older's inequality and of 
\eqref{revH}  yields
\begin{equation} \label{2017-2}
\bigg(\dashint_{ B_{\frac R2}(x)}\abs{\bfA(\nabla \bfu)}^{\min\set{p',2}}\dy \bigg)^\frac{1}{\min\set{p',2}}\lesssim\bigg(\dashint_{ B_{R}(x)}\abs{\nabla \bfu}\dy\bigg)^{p-1}+\bigg(\dashint_{ B_{R}(x)}\abs{\bfF-\mean{\bfF}_{ B_{R}(x)}}^{p'}\dy \bigg)^\frac{1}{p'}\,.
\end{equation}
 Finally, note that, by \eqref{tele4} 
\begin{align}\label{2017-3}
\sum_{i=0}^{k}\bigg(\dashint_{B_{\theta^i {\frac R2}}(x)}
\abs{\bfF-\mean{\bfF}_{B_{\theta^i {\frac R2}}}}^{p'}\dy\bigg)^\frac{1}{p'} & +\bigg(\dashint_{ B_{R}(x)}\abs{\bfF-\mean{\bfF}_{ B_{R}(x)}}^{p'}\dy \bigg)^\frac{1}{p'}\\
& \nonumber  \quad\lesssim
\int_r^{R} \bigg( \dashint_{B_\rho(x) }
 \abs{ \bfF-\mean{\bfF}_{B_\rho(x) }}^{p'} \dy \bigg)^\frac{1}{p'}  \, \frac{d\rho}{\rho}.
\end{align}
If $r \in (0, R/2]$, inequality \eqref{main2}  follows from  \eqref{ponc}, \eqref{2017-1},  \eqref{2017-2} and \eqref{2017-3}. Inequality \eqref{main2}  continues to hold  for $r \in (R/2, R]$, thanks to the Sobolev-Poincar\'e inequality.
\qed

\medskip
\par\noindent
\textbf{Proof of Theorem~\ref{mainpointwise}}.  Let $x \in \Omega$ and $R>0$ be such that $B_R(x) \subset \Omega$.
Lemma~\ref{lem:iter} and inequality \eqref{main2} yield
 \begin{align}\label{march2}
\bigg(\dashint_{B_{r}(x)}\abs{\bfu}dy\bigg)^{p-1} 
&\lesssim \bigg(\abs{\mean{\abs{\bfu}}_{B_{r}(x)}-\mean{\abs{\bfu}}_{B_{\frac R4}(x)}}
    + \dashint_{B_{\frac R4}(x)}\abs{\bfu}dy\bigg)^{p-1} \\
\nonumber &\lesssim
\bigg(\int_{r}^{\frac R4}\dashint_{B_{\rho}(x)}\abs{\bfu-\mean{\bfu}_{B_{\rho}(x)}}dy\frac{d\rho}{\rho}\bigg)^{p-1}
+ \bigg(\dashint_{B_{\frac R4}(x)}\abs{\bfu}dy\bigg)^{p-1}\\
\nonumber &\lesssim \bigg(\int _0^{\frac R4} \bigg(\int_\rho^{\frac R4}
\bigg(\dashint_{B_s(x)} \abs{\bfF-\mean{\bfF}_{B_s(x)}}^{p'} \dy
\bigg)^\frac{1}{p'}\frac{ds}{s}\bigg)^{\frac 1{p-1}} \,d\rho
\bigg)^{p-1}
\\
\nonumber & \quad +
  \bigg(\int _0^{\frac R4}\dashint_{B_{{\frac R4}}(x)}\abs{
\nabla\bfu}\dy\, dr\bigg)^{p-1}
%
%
\color{black}
 +
 \bigg(\dashint_{B_{\frac
R4}(x)}\abs{\bfu}dy\bigg)^{p-1},
\end{align}
for $\rho \in (0, R]$.
Next, observe that
\begin{equation}\label{march4}
\bigg(\int_{B_r(x)} \abs{\bfF-\mean{\bfF}_{B_r(x)}}^{p'} \dy
\bigg)^\frac{1}{p'} \leq 2 \bigg(\int_{B_r(x)}
\abs{\bfF-\mean{\bfF}_{B_\rho(x)}}^{p'} \dy \bigg)^\frac{1}{p'} \leq
2 \bigg(\int_{B_\rho(x)} \abs{\bfF-\mean{\bfF}_{B_\rho(x)}}^{p'} \dy
\bigg)^\frac{1}{p'}
\end{equation} 
 if $0<r\leq \rho<R$.  
As a consequence, the function
$$(0, R) \ni r \mapsto \bigg(\int_{B_r(x)} \abs{\bfF-\mean{\bfF}_{B_r(x)}}^{p'} \dy
\bigg)^\frac{1}{p'}$$ is quasi-increasing, with constant $k=2$, according to definition 
\eqref{quasi}. Hence, by Lemma \ref{hardy},   \color{black}
\begin{align}\label{march3}
\bigg(\int _0^{\frac R4} \bigg(\int_\rho^{\frac R4}
\bigg(\dashint_{B_s(x)} \abs{\bfF-\mean{\bfF}_{B_s(x)}}^{p'} \dy
\bigg)^\frac{1}{p'}\frac{ds}{s}\bigg)^{\frac 1{p-1}} \,d\rho
\bigg)^{p-1}
 \lesssim \bigg(
\int_0^{{\frac R2} \color{black}}\bigg( \dashint_{B_\rho(x)}
\abs{\bfF-\mean{\bfF}_{B_r(x)}}^{p'} \dy \bigg)^\frac{1}{p}d\rho
\bigg)^{p-1}.
\end{align}
From inequalities \eqref{march2} and \eqref{march3} one deduces that
\begin{align}\label{march5}
\dashint_{B_{\seb{r}}(x)}\abs{\bfu}dy   & \lesssim  \int_0^{\frac R2}\bigg(
\dashint_{B_\rho(x)} \abs{\bfF-\mean{\bfF}_{B_\rho(x)}}^{p'} \dy
\bigg)^\frac{1}{p}dr  \\ \nonumber & \quad \quad  + {R}\dashint_{B_{{\frac R4}}(x)}\abs{
\nabla\bfu}\dy  + 
\dashint_{B_{\frac R4}(x)}\abs{\bfu}\dy\,.
\end{align}
An inequality of Caccioppoli type tells us that 
\begin{align}
\label{cacc}
\bigg(\dashint_{B_{\frac R4}(x)}\abs{\nabla \bfu}^p\dy\bigg)^\frac1p\lesssim \frac 1 R\bigg(\dashint_{B_{\frac R2}(x)}\abs{\bfu-\mean{\bfu}_{B_{\frac R2}(x)}}^p\dy\bigg)^\frac1p + \bigg(\dashint_{B_{\frac R2}(x)}\abs{\bfF-\mean{\bfF}_{B_{\frac R2}(x)}}^{p'}\dy\bigg)^\frac1{p}\,,
\end{align}
see, for instance, \cite[Lemma 3.2]{DieKapSch11}. 
Inequality \eqref{cacc}, coupled with 
the {Sobolev}-\Poincare{} inequality, yields
\begin{align*}
\bigg(\dashint_{B_{\frac R4}(x)}\abs{\bfu-\mean{\bfu}_{B_{\frac R4}(x)}}^{qp}\dy\bigg)^\frac{1}{qp}\lesssim \bigg(\dashint_{B_{\frac R2}(x)}\abs{\bfu-\mean{\bfu}_{B_{\frac R2}(x)}}^p\dy\bigg)^\frac1p + R\bigg(\dashint_{B_{\frac R2}(x)}\abs{\bfF-\mean{\bfF}_{B_{\frac R2}(x)}}^{p'}\dy\bigg)^\frac1{p},
\end{align*}
for every
$q\in(1,\frac{n}{n-p}]$. Thus,
  \begin{align*}
\bigg(\dashint_{B_{\frac R4}(x)}\abs{\bfu}^{qp}\dy\bigg)^\frac{1}{qp}\lesssim \bigg(\dashint_{B_{\frac R2}(x)}\abs{\bfu}^p\dy\bigg)^\frac1p +   R\bigg(\dashint_{B_{\frac R2}(x)}\abs{\bfF-\mean{\bfF}_{B_{\frac R2}(x)}}^{p'}\dy\bigg)^\frac1{p},
\end{align*}
 whence, by and interpolation argument as in \cite[Remark 6.12]{Giu03}, one can deduce that
\begin{align}\label{revHu}
\bigg(\dashint_{B_{\frac R4}(x)}\abs{\bfu}^{p}\dy\bigg)^\frac{1}{p} \color{black} \lesssim  \dashint_{B_{\frac R2}(x)}\abs{\bfu}\dy +  R \bigg(\dashint_{B_{\frac R2}(x)}\abs{\bfF-\mean{\bfF}_{B_{\frac R2}(x)}}^{p'}\dy\bigg)^\frac1{p}.
\end{align}
Inequalities  \eqref{cacc} and \eqref{revHu} entail that
\begin{equation}\label{2017-8}
R\dashint_{B_{{\frac R4}}(x)}\abs{
\nabla\bfu}\dy \lesssim \dashint_{B_{\frac R2}(x)}\abs{\bfu}\dy + R\bigg(
\dashint_{B_{\frac R2}(x)} \abs{\bfF-\mean{\bfF}_{B_{\frac R2}(x)}}^{p'} \dy
\bigg)^\frac{1}{p}.
\end{equation}
Owing to equation   \eqref{march4} \color{black} 
\begin{equation}\label{2017-9}
R\bigg(
\dashint_{B_{\frac R2}(x)} \abs{\bfF-\mean{\bfF}_{B_{\frac R2}(x)}}^{p'} \dy \bigg)^\frac{1}{p}\lesssim \int_{\frac R2}^{R}\bigg(
\dashint_{B_\rho (x)} \abs{\bfF-\mean{\bfF}_{B_\rho (x)}}^{p'} \dy
\bigg)^\frac{1}{p}d\rho.
\end{equation}
Combining inequalities \eqref{march5}, \eqref{2017-8}, \eqref{2017-9}, and 
passing to the limit as $r \to 0^+$ yield inequality
\eqref{main1} for a.e. $x \in \Omega$, inasmuch as 
\begin{equation}\label{Wineq}
\int_0^{R}\bigg(
\dashint_{B_\rho(x)} \abs{\bfF-\mean{\bfF}_{B_\rho(x)}}^{p'} \dy
\bigg)^\frac{1}{p}d\rho \leq c \bfW_{\frac p{p+1}, p+1}^R(|\bfF|^{p'})(x) \quad \hbox{for $x\in \Omega$.}
\end{equation}
It remains to show that each point $x\in \Omega$ such that 
\begin{equation}\label{Leb}
\bfW_{\frac p{p+1}, p+1}^R(|\bfF|^{p'})(x) < \infty
\end{equation}
for some $R>0$ is a Lebesgue point of $\bfu$. Owing to inequality \eqref{Wineq}, assumption \eqref{Leb} implies that
\[
 \bigg(
\int_0^R\bigg( \dashint_{B_\rho(x)}
\abs{\bfF-\mean{\bfF}_{B_		\rho (x)}}^{p'} \dy \bigg)^\frac{1}{p}d\rho
\bigg)^{p-1}<\infty
\]
for some $R$. 
Let $0<s\leq r\leq R$.
By inequalities \eqref{tele1} and \eqref{main2},
\begin{align}\label{2017-10}
\abs{\mean{\bfu}_{B_r (x)}-\mean{\bfu}_{B_s(x)}} &
\lesssim \int_s^r\dashint_{B_{\rho}(x)}\abs{\bfu-\mean{\bfu}_{B_{\rho}(x)}}\,dy\,\frac{d\rho}{\rho}
\\ \nonumber
&\lesssim 
(r -s) \bigg[ \bigg(\int_0^R\bigg( \dashint_{B_\rho (x)}
\abs{\bfF-\mean{\bfF}_{B_\rho (x)}}^{p'} \dy
\bigg)^\frac{1}{p'}\frac{d\rho}{\rho}\bigg)^{\frac 1{p-1}} +  
\dashint_{B_{R}(x)}\abs{ \nabla\bfu}\dy\bigg].
\end{align}
Inequality \eqref{2017-10} ensures that the function $s \mapsto \mean{\bfu}_{B_s (x)}$ satisfies Cauchy's criterion, and hence converges to a finite limit $\bfu(x)$ as $s \to 0^+$. 
Next, by inequalities \eqref{main2} and \eqref{2017-10},
\begin{align}\label{2017-11}
\dashint_{B_r(x)}\abs{\bfu-\mean{\bfu}_{B_s (x)}}\, dy &\leq 
\dashint_{B_r(x)}\abs{\bfu-\mean{\bfu}_{B_r (x)}}\, dy
+\abs{\mean{\bfu}_{B_r (x)}-\mean{\bfu}_{B_s(x)}}
\\ \nonumber
&\lesssim
r \bigg(\int_0^{ R}
\bigg(\dashint_{B_\rho(x)} \abs{\bfF-\mean{\bfF}_{B_\rho(x)}}^{p'} \dy
\bigg)^\frac{1}{p'}\frac{d\rho}{\rho}\bigg)^{\frac 1{p-1}}
\\ \nonumber
&\quad + r  \bigg[\bigg(\int_0^R\bigg( \dashint_{B_\rho (x)}
\abs{\bfF-\mean{\bfF}_{B_\rho (x)}}^{p'} \dy
\bigg)^\frac{1}{p'}\frac{d\rho}{\rho}\bigg)^{\frac 1{p-1}} +  
\dashint_{B_{R}(x)}\abs{ \nabla\bfu}\dy\bigg].
\end{align}
Passing to the limit in \eqref{2017-11} first as $s \to 0^+$, and then as $r \to 0^+$ tells us that
$$\lim _{r \to 0^+}\dashint_{B_r(x)}\abs{\bfu-\bfu (x)} \dy
=0\,,$$
namely, $x$ is a Lebesgue point for $\bfu$. \qed

\section{Norm estimates}\label{appl}

The  Havin-Maz'ya-Wulff potential ${\bf W}_{\alpha , s}\bff$, and hence  ${\bf W}_{\alpha , s}^R\bff$ for every $R>0$, admit an estimate from
above, independent of $\bff$, in terms of the nonlinear potential
$\bfV_{\alpha , s}\bff$ given,
for $s \in  (1, \infty )$ and $\alpha \in (0,n)$, by
\begin{equation}\label{nonlinear}
{\bf V}_{\alpha , s}\bff (x) = {\bf I_\alpha} \big(({\bf I_\alpha}
|\bff|)^{\frac 1{s-1}}\big)(x) \qquad \quad \hbox{for $x \in \rn$.}
\end{equation}
Here, ${\bf I_\alpha}$ stands for the standard Riesz potential of order $\alpha$.
If the domain of $\bff$ is just  an open  set $\Omega\subset\rn$, then ${\bf W}_{\alpha , s}\bff$ and $\bfV_{\alpha , s}\bff$ are defined by continuing $\bff$ by $0$ outside $\Omega$.

The potentials ${\bf V}_{\alpha , s}$ extend the Riesz potentials, since, if $2 \alpha <
n $, then ${\bf V}_{\alpha , 2}\bff = c {\bf I}_{2\alpha }|\bff|$ for a
suitable constant $c=c(\alpha , n)$, and for every measurable
function $\bff$. They were introduced   by
V.P.Havin and V.G.Maz'ya \cite{HM}, and extensively investigated in the
framework of nonlinear capacity theory.
\\
By \cite[Theorem 3.1]{AM}, if $\alpha
s <n$, then there exists a constant $C=C(\alpha, s, n)$ such that
\begin{equation}\label{dis}
\bfW_{\alpha , s}\bff (x) \leq C \bfV_{\alpha , s}\bff (x) \qquad \quad
\hbox{for $x\in \rn$.}
\end{equation}
Incidentally, let us mention that a reverse pointwise inequality
 only holds for  $s$, $n$ and $\alpha$ in a suitable range -- see
\cite{HM}. However, a lower bound in integral form for $\bfW_{\alpha
, p}\bff $ in terms of $\bfV_{\alpha , s}\bff$ is  always available, as
proved in \cite{HW}.
\par As a consequence of inequalities \eqref{main1}  and  \eqref{dis},
local estimates for quasi-norms of solutions  to system \eqref{eq:sysA}, that are monotone with respect to pointwise ordering,   are reduced to boundedness properties of the nonlinear
potential $\bfV_{\frac{p}{p+1} , p+1}$ with respect to the same quasi-norms.
 We shall exhibit  the  estimates in question for  the Lorentz quasi-norms and the Orlicz norms. 
\\ Recall that, if $|\Omega|< \infty$, and either $q \in (1, \infty ]$, $\varrho \in (0, \infty]$, $\beta \in \mathbb R$, or $q=1$,  $\varrho \in (0, 1]$, $\beta \in [0, \infty)$, the Lorentz-Zygmund space $L^{q,\varrho}(\log L)^\beta (\Omega)$ is defined as the set of all measurable functions $\bff$ on $\Omega$ such that the quasi-norm given by
$$\|\bff \|_{L^{q,\varrho}(\log L)^\beta (\Omega)}= \|s^{\frac 1q - \frac 1{\varrho}}(1+\log (|\Omega|/s)^\beta|\bff|^*(s)\|_{L^\varrho (0, |\Omega|)}$$
is finite. Here, $|\bff|^*: [0, \infty) \to [0, \infty]$ denotes the decreasing--rearrangement of $|\bff|$. When $\beta =0$, the space $L^{q,\varrho}(\log L)^0 (\Omega)$ is a standard Lorentz space,   denoted by $L^{q,\varrho} (\Omega)$.  
In particular, 
$L^{q,q}(\log L)^{0}=L^{q,q}(\Omega) = L^q(\Omega)$ for $q \in [1, \infty]$. The notation $L^{q,\varrho}(\log L)^\beta _{\rm loc}(\Omega)$
stands for the set of functions which belong to $L^{q,\varrho}(\log L)^\beta (\Omega')$ for every open set $\Omega ' \subset \subset \Omega$.

\begin{theorem}\label{lorentz} Assume that $1<p<n$. Let $\bfu$ be a local weak solution to system \eqref{eq:sysA}.
\par\noindent
{\rm (i)} Let $0<\varrho \leq \infty$ and $1< q < \frac{n}{ p}$. If $\bfF
\in L^{qp', \varrho p'}_{\rm loc} (\Omega)$, then $\bfu \in
L^{\frac{qnp}{n-qp}, \varrho p}_{\rm loc} (\Omega)$.
\par\noindent
{\rm (ii)} Let $\varrho > \frac{1}{p}$.  If $\bfF \in L^{\frac n{p-1},
\varrho p'}_{\rm loc} (\Omega)$, then $\bfu \in L^{\infty, \varrho p}(\log
L)^{-1}_{\rm loc} (\Omega)$.
\par\noindent
{\rm (iii)} Let $0< \varrho \leq \frac{1}{p}$.  If $\bfF \in L^{\frac
n{p-1}, \varrho p'}_{\rm loc} (\Omega)$, then $\bfu \in L^{\infty}_{\rm loc}
(\Omega)$.
\end{theorem}

\begin{remark} {\rm Under the assumptions of Theorem \ref{lorentz}, Part (ii), one has that $\bfu \in \exp L^{(\varrho p)'}_{\rm loc} (\Omega)$, thanks to the  inclusion of the Lorentz-Zygmund space 
$L^{\infty, \varrho }(\log
L)^{-1} _{\rm loc}(\Omega)
$, with $\varrho >1$, into the Orlicz space of exponential type $\exp L^{\varrho '}_{\rm loc} (\Omega)$, whose definition is recalled below.
}
\end{remark}

A proof of Theorem \ref{lorentz} combines Theorem \ref{mainpointwise}  with the following characterization of 
boundedness properties of the potential $\bfV_{\alpha , s}\bff$ in  Lorentz type spaces from   \cite{Cianchinonlinear}.

\medskip
\par\noindent
{\bf Theorem A.} \, {\em
Let $s>1$, $n \geq 2$, and $0<\alpha < \frac ns$.  Let $\Omega$ be a
measurable subset of $\rn$.
\par\noindent
{\rm (i)} If $0<\varrho \leq \infty$ and $1< \sigma < \frac{n}{\alpha s}$,
then there exists a constant $C=C(\alpha , s, n, \sigma, \varrho )$
such that
\begin{equation}\label{10}
\|{\bf V}_{\alpha , s}\bff\|_{L^{\frac{\sigma n (s-1)}{n-\sigma \alpha
s}, \varrho (s-1)} (\Omega )} \leq C \|\bff\|_{L^{\sigma ,
\varrho}(\Omega )}^{\frac 1{s-1}}
\end{equation}
for every $\bff \in L^{\sigma , \varrho}(\Omega )$.
\par\noindent
{\rm (iii)} If $|\Omega |<\infty$ and $\varrho > \frac{1}{s-1}$,
 then there exists a
constant $C=C(\alpha ,  s, n,  \varrho , |\Omega |)$ such that
\begin{equation}\label{11bis}
\|{\bf V}_{\alpha , s}\bff\|_{L^{ \infty , \varrho (s-1)} (\log L)^{-1}
(\Omega )} \leq C \|\bff\|_{L^{\frac{n}{\alpha s} , \varrho}(\Omega
)}^{\frac 1{s-1}}
\end{equation}
for every $\bff \in L^{\frac{n}{\alpha p} , \varrho}(\Omega )$.
\par\noindent
(iv) If  $0< \varrho \leq \frac 1 {s-1}$, then there exists a
constant $C=C(\alpha , s, n)$ such that
\begin{equation}\label{12}
\|{\bf V}_{\alpha , s}\bff\|_{L^{ \infty} (\Omega )} \leq C
\|\bff\|_{L^{\frac{n}{\alpha s} , \varrho}(\Omega )}^{\frac 1{s-1}}
\end{equation}
for every $\bff \in L^{\frac{n}{\alpha s} , \varrho}(\Omega )$.}

\medskip
\par\noindent
{\bf Proof of Theorem \ref{lorentz}}. (i) Given any ball $B_R \subset \Omega$, we have, by
Theorem A, Part (i),
\begin{align}\label{lorentz1}
\big\|\bfV_{\frac p{p+1}, p+1}
\big(|\bfF|^{p'}\big)\big\|_{L^{\frac{qnp}{n-qp}, \varrho p}(B_R)} &
=
 \big\|\bfV_{\frac p{p+1}, p+1}
\big(|\bfF|^{p'}\big)\big\|_{L^{\frac{qnp}{n-\frac{qp}{p+1}(p+1)},
\varrho p}(B_R)} \\ \nonumber & \lesssim \big\||\bfF|^{p'} \big\|_{L^{q,
\varrho }(B_R)}^{\frac 1p} = \big\|\bfF  \big\|_{L^{qp', \varrho
p'}(B_R)}^{\frac 1{p-1}}.
\end{align}
The claim hence follows, via \eqref{main1} and \eqref{dis}.
\\
(ii)  Theorem A, Part (ii), ensures that
\begin{align}\label{lorentz2}
\big\|\bfV_{\frac p{p+1}, p+1}
\big(|\bfF|^{p'}\big)\big\|_{L^{\infty, \varrho p}(\log
L)^{-1}(B_R)}   \lesssim \big\||\bfF|^{p'} \big\|_{L^{\frac np, \varrho
}(B_R)}^{\frac 1p} = \big\|\bfF  \big\|_{L^{\frac n{p-1}, \varrho
p'}(B_R)}^{\frac 1{p-1}}.
\end{align}
The conclusion is now again a consequence of \eqref{main1} and \eqref{dis}.
\\
(iii)  By Theorem A, Part (iii),
\begin{align}\label{lorentz3}
\big\|\bfV_{\frac p{p+1}, p+1}
\big(|\bfF|^{p'}\big)\big\|_{L^{\infty}(B_R)}   \lesssim
\big\||\bfF|^{p'} \big\|_{L^{\frac np, \varrho }(B_R)}^{\frac 1p} =
\big\|\bfF  \big\|_{L^{\frac n{p-1}, \varrho p'}(B_R)}^{\frac
1{p-1}}.
\end{align}
A combination of inequalities  \eqref{main1} and \eqref{dis} completes the proof. \qed

\bigskip

Let us next focus on estimates involving norms in Orlicz spaces. Orlicz spaces are built upon Young functions, namely functions 
 $A: [0, \infty ) \to [0,
\infty ]$ that are convex (non trivial), left-continuous, and
vanish at $0$.
%
The Orlicz space $L^A(\Omega)$ consists of all measurable function $\bff$ on $\Omega$ for which the
 Luxemburg  norm 
\begin{equation*}
\|\bff\|_{L^A(\Omega)}= \inf \left\{ \lambda >0 :  \int_\Omega  A \left(
\frac{|\bff(x)|}{\lambda} \right) \dx \leq 1 \right\}
\end{equation*}
is finite. The local Orlicz space $L^A_{\rm loc}(\Omega)$ is defined accordingly. Observe that  $L^A (\Omega)=
L^q (\Omega)$ if $A(t)= t^q$ for some $q \in [1, \infty )$, and
$L^A (\Omega)= L^\infty (\Omega)$ if $A(t)= 0$ for $t\in [0,1]$, and $A(t)= \infty$ if $t \in (1, \infty)$. The Orlicz space associated with a Young function $A(t)$ which behaves like $t^p(\log t)^\alpha$ as $t \to \infty$, where either $p>1$ and $\alpha \in \mathbb R$, or $p=1$ and $\alpha \geq 0$, agrees with the Zygmund space denoted by $L^p(\log L)^\alpha (\Omega)$. The exponential type space $\exp L^\beta (\Omega)$, with $\beta >0$, is the Orlicz space built upon the Young function $A(t)=e^{t ^\beta} -1$. The double exponential space $\exp \exp L^{\beta}(\Omega)$ is defined analogously.
\par
Given two
Young functions $A$ and $B$, and $p \in (1,n)$,  define the
functions $E_{p}$ and $ F_{ p}$ from $[0, \infty )$ into $[0,
\infty]$ as
\begin{equation}\label{6}
E_{p}(t) = \bigg(\int _0^t \bigg(\frac {\tau ^{\frac 1{p} -
\frac{1}{n'}}}{A(\tau
)^{\frac{1}{n}}}\bigg)^{n'}d\tau\bigg)^{\frac{p}{n'}}
\end{equation}
and
\begin{equation}\label{7}
F_{p}(t) = \bigg(\int _0^t \frac {B(\tau )}{\tau ^{1 + \frac
n{n-p}}}d\tau\bigg)^{\frac{n- p}{n}}
\end{equation}
for $t\geq 0$.
\\ Our regularity result in Orlicz spaces reads as follows.
\begin{theorem}\label{orlicz}
Let $1 <p <n$,  and let $A$ and $B$ be
Young functions such that the functions   $E_p$ and $F_p$  defined as in \eqref{6} and \eqref{7} are finite-valued. Assume that 
there exist $\gamma
>0$ and $t_0 >0$ such that
\begin{equation}\label{8}
F_{p}\bigg(\frac{E_{ p}(t)}{\gamma } \bigg) \leq \gamma
\frac{A(t)}{t} \qquad \hbox{for $t>t_0$.}
\end{equation}
If   $|\bfF|^{p'} \in L^A_{\rm loc}(\Omega)$, and $\bfu$ is any local
weak solution  to system \eqref{eq:sysA},
then $|\bfu|^p \in L^B_{\rm loc}(\Omega)$.
%
%
\end{theorem}

\begin{remark}\label{convint}{\rm The assumption that the functions $E_p$ and $F_p$ are finite-valued, namely that the integrals on the right-hand sides of equations \eqref{6} and \eqref{7} are finite for $t\geq 0$, is immaterial in Theorem \ref{orlicz}. Indeed, $A$ and $B$ can be replaced, if necessary, by Young functions that are equivalent to $A$ and $B$ near infinity, and make the integrals in question converge. Indeed, such a replacement leaves the spaces $L^A_{\rm loc}(\Omega)$ and  $L^B_{\rm loc}(\Omega)$ unchanged.}
\end{remark}
\medskip

\begin{example}\label{exorlicz}
{\rm We specialize here Theorem \ref{orlicz} to the case when
$\bfF \in L^q(\log L)^\beta_{\rm
loc}(\Omega)$.
Assume that
$q>p'$ and $\beta \in \R$. An application of Theorem 
\ref{orlicz} tells us that, if  $\bfu$ is any local weak
solution  to system \eqref{eq:sysA}, then
%
\begin{equation*}
\bfu \in
\begin{cases}
L^{\frac {nq(p-1)}{n-q(p-1)}}(\log L)^{\frac{n\beta}{n-q(p-1)}}_{\rm
loc}(\Omega)
 &
 \hbox{if $p'<q <\frac n{p-1}$ and $\beta \in \R$,}
 \\
 \exp L^{\frac{n}{n-1 - \beta }}_{\rm loc}(\Omega)
 & \hbox{if $q=\frac n{p-1}$,  and $\beta < n-1$,} \\
 \exp \exp L^{{n'}}_{\rm  loc}(\Omega)
  & \hbox{if $q=\frac n{p-1}$,  and $\beta = n-1$}
  \\ L^\infty_{\rm loc}(\Omega) & \hbox{if either $q >\frac
n{p-1}$ and $\beta \in \R$, or $q=\frac n{p-1}$  and $\beta
> n -1$}
  .
\end{cases}
\end{equation*}
}
\end{example}

Theorem \ref{orlicz} is a consequence of  Theorem \ref{mainpointwise} and of (a special case of) a characterization of boundedness properties of the potential 
 $({\bf V}_{\alpha , s})^{p-1}$ in Orlicz spaces established in 
\cite{Cianchinonlinear}, and recalled in Theorem B below. Its statement involves the functions  $E_{\alpha , s}$ and $ F_{\alpha , s}$ associated with the Young functions $A$ and $B$ as follows.
Let $0< \alpha <\min\{\frac ns, \frac n{s'}\}$.
Define  the functions $E_{\alpha , s}$ and $ F_{\alpha , s} : [0,
\infty ) \to [0, \infty]$ as
\begin{equation}\label{6V}
E_{\alpha , s}(t) =
\bigg(\int _0^t \bigg(\frac {\tau ^{\frac 1{s-1} -1
+\frac{\alpha s'}{n}}}{A(\tau )^{\frac{\alpha s'}{n}}}\bigg)^{\frac
n{n- \alpha s'}}d\tau\bigg)^{s-1-\frac{\alpha s}{n}}\,,
\end{equation}
and
\begin{equation}\label{7V}
F_{\alpha , s}(t) = \bigg(\int _0^t \frac {B(\tau )}{\tau ^{1 +
\frac n{n-\alpha s}}}d\tau\bigg)^{\frac{n-\alpha s}{n}}
\end{equation}
for $t\geq 0$.

\medskip
\par\noindent
{\bf Theorem B.} \, {\em
Let $s>1$, $n \geq 2$, and $0<\alpha < \min\{\frac ns, \frac n{s'}\}$.  Let $\Omega$ be a
measurable subset of $\rn$ with $|\Omega |< \infty$. Let $A$
and $B$ be Young functions such that the  functions  $E_{\alpha , s}$ and $F_{\alpha , s}$,  defined as in 
\eqref{6V} and \eqref{7V}, are  finite-valued. Assume that there exist $\gamma
>0$ and $t_0 >0$ such that
\begin{equation}\label{8V}
F_{\alpha , s}\bigg(\frac{E_{\alpha , s}(t)}{\gamma } \bigg) \leq
\gamma \frac{A(t)}{t} \qquad \hbox{for $t>t_0$.}
\end{equation}
Then there exists a constant $C=C(\alpha ,  s, n, \gamma , A)$ such
that
\begin{equation}\label{9V}
\|({\bf V}_{\alpha , s}\bff)^{s-1}\|_{L^B(\Omega )} \leq C
\|\bff\|_{L^A(\Omega )}
\end{equation}
for every $\bff \in L^A(\Omega )$.
}

\begin{remark}\label{convintV} {\rm 
 The assumption on the finiteness of the functions $E_{\alpha, s}$ and $F_{\alpha, s}$ in Theorem B is irrelevant, by a reason analogous to that explained in Remark \ref{convint}.}
\end{remark}

\medskip
\par

In the remaining part of this section, we present  a few applications of Theorem \ref{mainosc} to the
regularity of solutions to \eqref{eq:sysA} in spaces of Campanato
type, and, as a consequence, in $\setBMO$ and in spaces of uniformly continuous functions.
\\
Given an exponent $q \geq 1$ and a continuous function  $\omega : (0, \infty ) \to (0, \infty)$, the Campanato type space $\mathcal L ^{\omega (\cdot ), q}(\Omega)$   is the space of those  functions $\bff \in L^q_{\rm loc}(\Omega )$ for which
the semi-norm
\begin{align}\label{campnorm}
  \norm{\bff}_{\mathcal L ^{\omega (\cdot ), q}(\Omega)}= \sup _{B_r \subset \Omega}
  \frac 1{\omega (r)}\bigg(\dashint_{B_r} |\bff - \mean{\bff}_{B_r}|^q\,
  dx \bigg)^{\frac 1q}
%
\end{align}
is finite.  The space $\mathcal L ^{\omega (\cdot ) ,q} _{\rm loc}(\Omega)$
is defined accordingly, as the set of all functions $\bff$ such that
$ \norm{\bff}_{\mathcal L ^{\omega (\cdot ), q}(\Omega ')} < \infty$ for every
open set $\Omega ' \subset \subset \Omega$. 
%
%
  Observe that $\mathcal L^{\omega (\cdot), q}_{\rm loc}(\Omega)$ depends only on the behavior of $\omega$ in a right neighborhood of $0$. The same is true for  $\mathcal L^{\omega (\cdot), q}(\Omega)$, if $\Omega$ is bounded.
\\
When $q=1$, we 
simply denote  $\mathcal L ^{\omega (\cdot ), 1}  (\Omega)$ by  $\mathcal L
^{\omega (\cdot )} (\Omega)$, and similarly
for local spaces.
 Let us notice that, 
 if $\omega$ is non-decreasing, then $$\mathcal L ^{\omega (\cdot ), q}_{\rm loc}(\Omega)= \mathcal L
^{\omega (\cdot )} _{\rm loc}(\Omega)$$ for every $q \geq 1$. This is an easy consequence of the John-Nirenberg theorem, asserting that,
if $\omega (r) =1$ and $q \geq 1$, then
\begin{align}\label{bmonorm}\mathcal L ^{\omega (\cdot),q}_{\rm loc}(\Omega) = \setBMO _{\rm loc}(\Omega)\,,
\end{align}
the space of functions of (locally) bounded mean oscillation. Another standard 
 choice is 
 \begin{equation}\label{omega}
  \omega (r) = r^\beta \quad \hbox{for $r \geq 0$,}
 \end{equation}
 for some   $\beta
\leq 1$\color{black}, in which case we also make  use of the notation $\mathcal L^{\beta }(\Omega)$ instead of $\mathcal L^{\omega (\cdot) }(\Omega)$. 
\\ We next denote by
$C^{\omega (\cdot)}(\Omega)$ the space of those functions $\bff$ in $\Omega$ such that
$$
\sup_{
\begin{tiny}
 \begin{array}{c}{
    x, y \in \Omega} \\
x\neq y
 \end{array}
  \end{tiny}
}\frac{|\bff (x) - \bff (y)|}{\omega (|x-y|)} < \infty\,.$$ 
The notation $C^{\omega (\cdot)}_{\rm loc}(\Omega)$ is employed accordingly.
Thus,  if $\lim _{r \to 0^+} \omega (r) =0$, then
$C^{\omega (\cdot)}(\Omega)$ is the space
of uniformly continuous functions, with modulus of continuity not
exceeding $\omega$. In particular, if $\omega$ has the form \eqref{omega} for some $\beta \in (0,1]$, then $C^{\omega (\cdot)}(\Omega)$ agrees with the space of H\"older continuous functions in $\Omega$, with exponent $\beta$. The latter is denoted by $C^\beta (\Omega)$  when $\beta \in (0,1)$, and by    ${\rm Lip}(\Omega)$ when $\beta =1$. 
 \\
 More generally, it is easily verified that 
\begin{equation}\label{august106}
C^{\omega (\cdot)}_{\rm loc}(\Omega) \subset  \mathcal L ^{\omega  (\cdot)}_{\rm loc}(\Omega)\,.
\end{equation}
 However, the reverse inclusion in \eqref{august106} need not hold for an arbitrary
$\omega$, as shown, for example, by equation \eqref{bmonorm}.
Results from \cite{Spa65} ensure that, if $\omega $  decays at $0$ so fast that
\begin{equation}\label{dini}
\int _0 \frac{\omega (r)}r\, dr < \infty,
\end{equation}
then
\begin{equation}\label{spanne1}
\mathcal L ^{\omega  (\cdot)}_{\rm loc}(\Omega) \subset  C^{\varpi  (\cdot)}_{\rm loc}(\Omega)\,,
\end{equation}
where the function $\varpi : [0,
\infty) \to [0, \infty )$ is defined as
\begin{equation}\label{spanne2}
\varpi (r) = \int _0^r \frac{\omega (\rho)}\rho\, d\rho \qquad
\hbox{for $r \geq 0$.}
\end{equation}
For instance, if $\omega$ is given by \eqref{omega}, then \eqref{spanne1} recovers Campanato's representation theorem asserting that
\begin{equation}\label{holdernorm}
\mathcal L^\beta _{\rm loc}(\Omega) = \begin{cases}
C^\beta _{\rm loc}(\Omega) & \hbox{if $\beta \in (0,1)$}
\\ 
{\rm Lip}_{\rm loc}(\Omega) & \hbox{if $\beta=1$.}
\end{cases}
\end{equation}
On the other hand, if \eqref{dini} fails, and the function $r
\mapsto \tfrac{\omega (r)}r$ is non-increasing near $0$, then the
space $\mathcal L ^{\omega (\cdot)}_{\rm loc}(\Omega)$ is not even contained in
$L^\infty _{\rm loc} (\Omega)$.
\\
The  Morrey type space $\mathcal M^{\omega (\cdot) , q}(\Omega)$ is
 defined as the space of those functions
 $\bff \in L^q_{\rm loc}(\Omega )$ such that the norm
\begin{align}\label{morreynorm}
  \norm{\bff}_{\mathcal M ^{\omega (\cdot), q}(\Omega)}= \sup _{B_r \subset \Omega}
  \frac 1{\omega (r)} \bigg(\int_{B_r} |\bff |^q\,
  dx\bigg)^{\frac 1q}
%
\end{align}
is finite. The local Morrey space 
$\mathcal M^{\omega (\cdot), q}_{\rm loc}(\Omega)$ is defined with the usual modification.
When $\omega$ has the form \eqref{omega}, we also make use of the alternate notation $\mathcal M^{\beta, q}(\Omega)$ instead of $\mathcal M^{\omega (\cdot), q}(\Omega)$.

\medskip

The regularity of solutions $\bfu$ to system \eqref{eq:sysA} in Campanato spaces, for data $\bfF$ in Morrey spaces, is described in the following result.

\begin{theorem}\label{campestimate}
Let $\omega : (0, \infty) \to (0, \infty)$ be a continuous
function, and let  $\mu : (0, \infty) \to (0, \infty)$ be a
function given by
$$\mu (r) = r \bigg(\int _r^ 1 \frac{\omega (\varrho)}{\varrho
^{\frac{n}{p'}+1}} \, d\varrho\bigg)^{\frac 1{p-1}} $$
%
%
for $r$ near $0$. If $\bfF \in \mathcal M^{\omega (\cdot), p'}_{loc}(\Omega)$, and $\bfu$ is any
local weak solution  to system
\eqref{eq:sysA},  then $\bfu \in \mathcal L
^{\mu (\cdot)}_{loc}(\Omega)$.
\end{theorem}

\par\noindent
{\bf Proof.}
Inequality \eqref{main2} ensures that there exists a constant $C$ such that
\begin{align}\label{may1}
 \frac 1{\mu (r)} \dashint_{B_r(x)}\abs{\bfu(y)-\mean{\bfu}_{B_r (x)}}\, dy
&\leq \frac {Cr }{\mu (r)} \bigg(\int_r^R\bigg( \dashint_{B_\rho(x)}
\abs{\bfF-\mean{\bfF}_{B_\rho(x)}}^{p'} \dy
\bigg)^\frac{1}{p'}\,\frac{d\rho}{\rho}\bigg)^{\frac 1{p-1}}  \\ \nonumber & + \frac {Cr }{\mu (r)} 
\dashint_{B_{R}(x)}\abs{ \nabla\bfu}\dy
\end{align}
for every $x \in \Omega$ and $R>0$ such that $B_R(x) \subset \Omega$, and $r \in (0, R]$.  On the other hand,
\begin{align}\label{may2}
\frac {r }{\mu (r)} \bigg(\int_r^R& \bigg( \dashint_{B_\rho (x)}
\abs{\bfF-\mean{\bfF}_{B_\rho(x)}}^{p'} \dy
\bigg)^\frac{1}{p'}\frac{d\rho}{\rho}\bigg)^{\frac 1{p-1}} \\ \nonumber & \leq  
\frac {Cr }{\mu (r)} \bigg(\int_r^R\bigg( \int_{B_\rho(x)}
\abs{\bfF}^{p'} \dy
\bigg)^\frac{1}{p'}\frac{d\rho}{\rho ^{\frac n{p'}+1}}\bigg)^{\frac 1{p-1}} 
\\  \nonumber & 
\leq 
\frac {Cr }{\mu (r)} \bigg(\int_r^R 
 \frac{\omega (\rho)}{\rho^{\frac n{p'}+1}}\, d\rho \bigg)^{\frac 1{p-1}} \bigg(\sup _{\rho\in (0, R)} \frac 1{\omega (\rho)}\bigg(\int_{B_\rho(x)}\abs{\bfF}^{p'} \dy\bigg)^\frac{1}{p'}
\bigg)^\frac{1}{p-1} 
\end{align}
for some constant $C$ and for every $r \in (0, R]$. Hence, the conclusion follows. \qed

\begin{example}\label{expower}{\rm  
Assume that $\bfF \in \mathcal M^{\beta, p'}_{\rm loc}(\Omega)$ for some $ \beta \geq (n-p)/p'$. An application of Theorem \ref{campestimate} tells us the following. If 
$\beta =(n-p)/p'$, then 
$$\bfu \in \setBMO_{\rm loc}(\Omega);$$
if $(n-p)/p' < \beta  < n/p'$, then $$\bfu \in C^{1- (\frac np - \frac \beta{p-1})}_{\rm loc}(\Omega);$$
if $\beta = n/p'$, then $$\bfu \in C^{\nu (\cdot)}_{\rm loc}(\Omega)$$
with $\nu (r) \approx r(\log 1/r)^{\frac 1{p-1}}$ for $r$ near  $0$; if $\beta >n/p'$, then $$\bfu \in {\rm Lip}_{\rm loc}(\Omega).$$ \color{black}
}
\end{example}

\begin{corollary}\label{holder}
Assume that   $\bfF
\in L^q_{\rm loc}(\Omega)$ for some $q > \max \{p', \frac n{p-1}\}$. If $\bfu$ is any
local weak solution  to system
\eqref{eq:sysA},  
then $\bfu \in C^{1-\frac n{q(p-1)}}_{\rm loc}(\Omega)$.
\end{corollary}
\par\noindent
{\bf Proof.} H\"older's inequality ensures that, if $\bfF
\in L^q_{\rm loc}(\Omega)$, then $\bfF \in \mathcal M^{\frac {n(q-p')}{qp'}, p'}_{\rm loc}(\Omega)$. Theorem \ref{campestimate} then entails that $u \in \mathcal L ^{1-\frac n{q(p-1)}}_{\rm loc}(\Omega)$. Hence, the conclusion follows, owing to equation \ref{holdernorm}.
\qed

\bigskip
Our last result concerns the Lipschitz continuity of weak solutions to system \eqref{eq:sysA}.

\begin{theorem}\label{Lip} Let $\omega : (0, \infty) \to (0, \infty)$ be a continuous function fulfilling condition \eqref{dini}.
If $\bfF \in \mathcal L^{\omega (\cdot), p'} _{\rm loc}(\Omega)$, and $\bfu$ is any
local weak solution  to system
\eqref{eq:sysA},  then $\bfu \in {\rm Lip}_{\rm loc}(\Omega)$.
\end{theorem}

\smallskip
\par\noindent
{\bf Proof.} By inequality \eqref{mainosc}, there exists a constant $C$ such that
\begin{align}\label{Lip1}
\frac 1r \dashint _{B_r(x)}\abs{\bfu -\mean{\bfu}_{B_r (x)}} \dy & \leq C
\bigg(\int_r^R\bigg( \dashint_{B_\rho(x)}
\abs{\bfF-\mean{\bfF}_{B_\rho(x)}}^{p'} \dy
\bigg)^\frac{1}{p'}\frac{d\rho}{\rho}\bigg)^{\frac 1{p-1}}   + C
\dashint_{B_{R}(x)}\abs{ \nabla\bfu}\dy
\\ \nonumber & 
\leq C \|\bfF\|_{\mathcal L^{\omega (\cdot), p'} _{\rm loc}(B_R(x))} \bigg(\int _0^R \frac{\omega (\rho)}{\rho}\, d\rho\bigg)^{\frac 1{p-1}}
+ C
\dashint_{B_{R}(x)}\abs{ \nabla\bfu}\dy\,
\end{align}
for every $x \in \Omega$ and $R>0$ such that $B_R(x) \subset \Omega$, and $r \in (0, R]$. Thus, $\bfu \in \mathcal L^{\omega (\cdot)} _{\rm loc}(\Omega)$, with $\omega (r) =r$, whence $\bfu \in {\rm Lip}_{\rm loc}(\Omega)$, by \eqref{holdernorm}.
\qed

\bigskip\par\noindent
{\bf Acknowledgments}. \,
This research was partly supported by: Research project of MIUR
(Italian Ministry of Education, University and Research) Prin 2015,
n.2015HY8JCC,  \lq\lq Partial differential equations and related analytic-geometric inequalities";  GNAMPA of the Italian INdAM (National Institute of
High Mathematics);  PIA 2016-2017 of the University of Florence (Italy);
Project LL1202  financed by the Ministry of Education, Youth and Sports, Czech Republic.

\end{document}